\documentclass[lettersize,journal]{IEEEtran}
\usepackage{amsmath,amsfonts}
\usepackage{algorithmic}
\usepackage{algorithm}
\usepackage{array}
\usepackage{multirow}
\usepackage[caption=false,font=normalsize,labelfont=sf,textfont=sf]{subfig}
\usepackage{textcomp}
\usepackage{stfloats}
\usepackage{amssymb}
\usepackage{pifont}
\newcommand{\xmark}{\ding{55}}
\usepackage{bm}
\usepackage{url}
\usepackage{amsthm}
\usepackage{verbatim}
\usepackage{graphicx}
 \usepackage{booktabs}
\usepackage{cite}

\newtheorem{thm}{Theorem}

\newtheorem{defn}{Definition}

\usepackage{color}
\usepackage{comment}
\usepackage{booktabs}
\usepackage{caption}
\usepackage{subcaption}
\newtheorem{ex}{Example}
\usepackage{mathrsfs} 

\newtheorem{assum}{Assumption}
\newtheorem{coro}{Corollary}
\newtheorem{prf}{Proof}
\newcommand{\norm}[1]{\|#1\|}

\usepackage{mathrsfs}
\begin{document}

\title{Designing Robust Software Sensors for Nonlinear Systems via Neural Networks and Adaptive Sliding Mode Control
}
\author{%
  Ayoub~Farkane, 
  Mohamed~Boutayeb, 
  Mustapha~Oudani, 
  and~Mounir~Ghogho,,~\IEEEmembership{Fellow,~IEEE}%
  \thanks{A. Farkane is with TICLab, International University of Rabat, Rabat 11100, Morocco; \emph{and} with CNRS–CRAN–7039, University of Lorraine, 54000 Nancy, France (e-mail: ayoub.farkane@uir.ac.ma).}%
  \thanks{M. Boutayeb is with CNRS–CRAN–7039, University of Lorraine, 54000 Nancy, France; \emph{and} with INRIA Nancy–LARSEN, 54600 Villers-lès-Nancy, France.}%
  \thanks{M. Oudani is with TICLab, International University of Rabat, Rabat 11100, Morocco.}%
  \thanks{M. Ghogho is with the Mohammed VI Polytechnic University, College of Computing, 43150 Benguerir, Morocco; \emph{and} with   Faculty of Engineering, University of Leeds, LS2 9JT Leeds, United Kingdom(e-mail: mounir.ghogho@um6p.ma).}%
}


\IEEEpubid{}

\maketitle

\begin{abstract}
Accurate knowledge of the state variables in a dynamical system is critical for effective control, diagnosis, and supervision, especially when direct measurements of all states are infeasible. This paper presents a novel approach to designing software sensors for nonlinear dynamical systems expressed in their most general form. Unlike traditional model-based observers that rely on explicit transformations or linearization, the proposed framework integrates neural networks with adaptive Sliding Mode Control (SMC) to design a robust state observer under a less restrictive set of conditions. The learning process is driven by available sensor measurements, which are used to correct the observer's state estimate. The training methodology leverages the system's governing equations as a physics-based constraint, enabling observer synthesis without access to ground-truth state trajectories. By employing a time-varying gain matrix dynamically adjusted by the neural network, the observer adapts in real-time to system changes, ensuring robustness against noise, external disturbances, and variations in system dynamics. Furthermore, we provide sufficient conditions to guarantee estimation error convergence, establishing a theoretical foundation for the observer's reliability. The methodology's effectiveness is validated through simulations on challenging examples, including systems with non-differentiable dynamics and varying observability conditions. These examples, which are often problematic for conventional techniques, serve to demonstrate the robustness and broad applicability of our approach. The results show rapid convergence and high accuracy, underscoring the method's potential for addressing complex state estimation challenges in real-world applications.
\end{abstract}


\section{Introduction}
\IEEEPARstart{S}{t}ate estimation is a pivotal component  for the  control and supervision of dynamical systems, particularly when direct measurement of all state variables is infeasible due to cost, sensor limitations, or physical constraints. In nonlinear systems, estimating unmeasured states becomes more complex as state variables exhibit intricate, non-linear behavior. This complexity necessitates advanced observer designs that ensure accurate state estimation, especially under dynamic and uncertain conditions.

Traditional model-based observers, such as the Extended Kalman Filter (EKF), often rely on assumptions of linearity or transformations that may not hold in more complex systems. These methods can struggle to achieve accurate state tracking in the presence of system nonlinearity, noise, or external disturbances. Consequently, there is a growing need for adaptive approaches to overcome these limitations and provide robust state estimation.

Recent research \cite{buisson2023towards,peralez2021deep,ramos2020numerical,niazi2023learning, miao2023learning,alvarez2024nonlinear} have focused on leveraging the capabilities of deep neural networks to approximate nonlinear functions and apply them to the design and estimation of nonlinear dynamic systems. Many works employ neural networks to learn transformation maps and linearize the system. While these approaches can be effective for many systems, they may struggle with complex problems \cite{andrieu2021remarks,bernard2022observer}. This raises a crucial research question: can we develop a neural network architecture capable of directly learning complex nonlinear dynamics without explicit transformation or linearization?

This paper proposes a novel approach for state estimation in nonlinear dynamical systems that integrates neural networks with adaptive Sliding Mode Control (SMC). This approach utilizes a time-varying gain matrix that adapts to real-time changes in system dynamics, thereby enhancing estimation accuracy and stability. The neural network component dynamically learns optimal gain adjustments based only on the available sensor measurements or measured information, ensuring that state estimation is achieved without requiring full knowledge of the system dynamics. This eliminates the need for system linearization 
and makes the method widely applicable to complex and uncertain systems. The contributions of this work are  as follows:
(1) Novel Integration of Neural Networks with Adaptive SMC: We propose a new framework that leverages neural networks to dynamically learn and adjust time-varying gain matrices, ensuring accurate state estimation without requiring system linearization.
(2) Enhanced Robustness and Stability: The adaptive SMC mechanism is designed to handle dynamic system changes effectively, providing robustness against noise and external disturbances.
(3) Real-time Adaptability: Our approach enables real-time adaptation to the varying dynamics of the system, making it suitable for practical applications that require online learning and rapid convergence.
(4) Comprehensive Evaluation: We rigorously test the proposed method against state-of-the-art learning-based observers, under both noise-free and noisy conditions.

The structure of this paper is as follows: Section II reviews related works on state observers and adaptive control techniques. Section III introduces the mathematical formulation and challenges associated with nonlinear state estimation. Section IV presents the proposed neural network-based observer and its integration with adaptive SMC. Section V details the convergence analysis of the approach, while Section VI demonstrates its effectiveness through extensive applications and simulations. Section VII concludes the paper with key insights and potential directions for future research.
\section{Related works}
The design of state observers for nonlinear systems is a foundational problem in control theory, with a research landscape that has evolved from classical model-based methods to modern learning-based and hybrid frameworks.

\subsection{A Classical Model-Based and Adaptive Observers}

Traditional approaches are typified by the Extended Kalman Filter (EKF), which relies on local linearization of the system dynamics. While widely used, the EKF's performance is fundamentally limited by model accuracy and the validity of the linearization assumption, especially for highly nonlinear systems \cite{einicke1999robust,boutayeb1997convergence}. To handle parametric uncertainty, classical adaptive observers were developed, which can simultaneously estimate states and unknown parameters. However, many of these approaches require linear parameterizations and strict conditions like persistence of excitation. More advanced non-neural adaptive schemes, such as the adaptive sliding-mode observer \cite{rios2018adaptive}, have been proposed to achieve faster parameter convergence using nonlinear adaptation laws, but the challenge of handling unstructured, non-parametric uncertainties remains.

\subsection{Learning-Based Observers}

To overcome the dependency on precise models, machine learning has emerged as a powerful tool for observer design. These efforts have largely progressed along two parallel paths: learning a coordinate transformation to a simpler system, or directly learning the system dynamics or observer gains.
\begin{enumerate}
    \item  Learning Transformation Maps to a Linear Domain: A prominent strategy is to learn a transformation that maps the nonlinear system into a more tractable, often linear, domain. The Kazantzis-Kravaris-Luenberger (KKL) observer framework exemplifies this, where neural networks are trained to find an injective immersion into a linear system\cite{ramos2020numerical, buisson2023towards,peralez2021deep}. This concept has been significantly advanced recently. For instance, \cite{miao2023learning}  use Neural ODEs to design KKL observers and critically analyze the trade-off between convergence speed and robustness. \cite{marani2025deep} employ a deep-learning KKL chain observer to handle challenging time-varying output delays. Furthermore, Physics-Informed Neural Networks (PINNs) have been used by \cite{alvarez2024nonlinear} to numerically solve the functional equations inherent in this exact linearization framework. While elegant, the success of these methods hinges on the existence and, crucially, the learnability of such a complex transformation \cite{andrieu2021remarks,bernard2022observer}.
\item Directly Learning System Dynamics or Observer Gains: An alternative path is to use neural networks to directly approximate the unknown parts of the system dynamics or the observer gain itself. Methods employing Echo State Networks (ESNs) \cite{goswami2021data,maass2004computational,grigoryeva2018echo} or Neural ODEs \cite{miao2023learning} fall into this category. \cite{zhuang2023neural}  propose a Neural Network Adaptive Observer (NNAO) using RBFs to explicitly handle systems with partially or completely unknown dynamics subject to practical sampling and delay constraints. In a highly sophisticated approach, \cite{marani2024unsupervised} leverage PINNs to learn the observer's correction term by directly enforcing the stringent conditions of contraction analysis within the network's loss function, creating a robust, unsupervised design.
\end{enumerate}
\subsection{ The Need for Integrated, Robust, and Adaptive Architectures}

While the learning-based methods above offer remarkable flexibility, they often lack formally integrated mechanisms for robustness against external disturbances, sensor noise, and unmodeled dynamics. This has motivated a third, and arguably most promising, direction: the development of hybrid observers that combine the expressive power of machine learning with the proven stability and robustness guarantees of established control and filtering theories.

One class of hybrid observers fuses NNs with established filtering techniques. For example, \cite{de2024hybrid} integrate a PINN with an adaptive Unscented Kalman Filter (UKF), using the PINN to model the nonlinear dynamics and the UKF to handle uncertainty and noise systematically.

Another, highly active, class of hybrid observers integrates neural networks with the robust framework of Sliding Mode Control (SMC). The core idea is to use the NN to learn the unknown dynamics, thereby enabling the design of a high-performance SMC/SMO that guarantees stability. Recent works have explored this synergy from various angles. \cite{zhong2024fuzzy} combine a type-2 fuzzy neural network with adaptive SMC for uncertain descriptor systems. \cite{taimoor2024neural} and \cite{yang2025neural} also propose SMC-based observers, with the former using a fuzzy SMC as the learning algorithm and the latter focusing on event-triggered control for networked systems. The same methodology has been successfully applied to control, as demonstrated by \cite{vacchini2023design},  in which deep NNs are used to enable an Integral Sliding Mode Control (ISMC) scheme for systems with fully unknown dynamics.

These works collectively demonstrate a clear and powerful trend: leveraging NNs to learn the "what" (the system dynamics) so that robust control theory can provide the "how" ( a stability guarantee). However, many of these approaches either rely on specific structures like fuzzy logic or are tailored to particular problems like fault detection or networked systems. This highlights the need for a more generalized and versatile framework. Our work addresses this gap by proposing a nearly universal approach that integrates a neural network with adaptive SMC to learn a time-varying observer gain. This design does not rely on explicit system linearization or restrictive assumptions, providing a robust, accurate, and highly adaptable solution for a wide class of nonlinear dynamical systems.
\section{Governing system}
Let us consider the following nonlinear discrete-time dynamic system:
\begin{equation}
\left\{
\begin{aligned}
x_{k+1} &= f(x_k) + Bu_k, \\
y_k &= h(x_k), \\
x(0) &= x_0
\end{aligned}
\right. \label{e1}
\end{equation}
where $x_k \in \mathbb{R}^n$ is the system state at time step $k$, $f$ is a nonlinear function representing the system dynamics, $y_k \in \mathbb{R}^m$ is the measured output, $B \in \mathbb{R}^{n \times p}$ is a known input matrix,  $u_k \in \mathbb{R}^p$ is the control input, and  $x_0$ is the initial state.

The existence and uniqueness of system solutions \eqref{e1} are guaranteed for any given initial state $x_0$ by the Cauchy-Lipschitz theorem under the assumption of local Lipschitz continuity for the function $f$.

In many practical scenarios, direct measurement of all system states is infeasible due to cost, sensor limitations, or physical constraints. To address this, a state observer is employed to estimate the unmeasured or unmeasurable internal states based on available inputs and outputs. The significance of state observers is particularly evident in nonlinear systems, where certain states may exhibit complex, nonlinear behaviors that are challenging to measure directly.

\section{Proposed Model}\label{mo}
A conventional nonlinear state observer for system \eqref{e1} can be formulated as follows:

\begin{equation}
\left\{\begin{aligned}
\hat{x}_{k+1} &= f(\hat{x}_k) + L_k (y_k - \hat{y}_k) + Bu_k \\
\hat{y}_k& = h(\hat{x}_k) \\
\hat{x}(0)&=\hat{x}_0
\end{aligned}
\right.
\label{e2}
\end{equation}
where $\hat{x}_k \in \mathbb{R}^n$ denotes the estimated state, $\hat{y}_k \in \mathbb{R}^m$ represents the estimated output, and $L_k \in \mathbb{R}^{n \times m}$ is the observer gain matrix.

The primary objective of the nonlinear state observer in \eqref{e2} is to determine an optimal gain matrix $L_k$ that ensures the convergence of the estimated state $\hat{x}_k$ to the true state $x_k$, i.e., $\lim_{k \rightarrow \infty} \| x_k - \hat{x}_k \| = 0$. Employing a time-varying gain matrix offers significant advantages over a constant gain. This adaptive approach enables the observer to effectively handle time-varying system dynamics, leading to improved state estimation accuracy under changing operating conditions.

To enhance robustness, we introduce an adaptive SMC term, denoted as $\nu_k$, to the observer \eqref{e2}. The modified observer structure is then given by:
\begin{equation}
\left\{\begin{aligned}
\hat{x}_{k+1} &= f(\hat{x}_k) + L_k (y_k - \hat{y}_k) + Bu_k + \nu_k \\
\hat{y}_k& = h(\hat{x}_k) \\
\hat{x}(0)&=\hat{x}_0
\end{aligned}
\right.
\label{e3}
\end{equation}
To address the challenge of determining the optimal time-varying gain matrix $L_k$, we leverage a neural network to approximate its values. The network's parameters are continuously updated based on real-time measurements ($y_k$) and control inputs ($u_k$). This adaptation process is guided by a specifically designed loss function. The loss function incorporates the nonlinear system dynamics and aims to minimize the estimation error (i.e., the difference between the actual output $y$ and the estimated output $\hat{y}$). By minimizing this loss, the neural network is effectively trained to generate optimal gain matrices that enhance the overall performance of the state observer.
 \begin{figure*}[!t]
     \centering     \includegraphics[width=\textwidth]{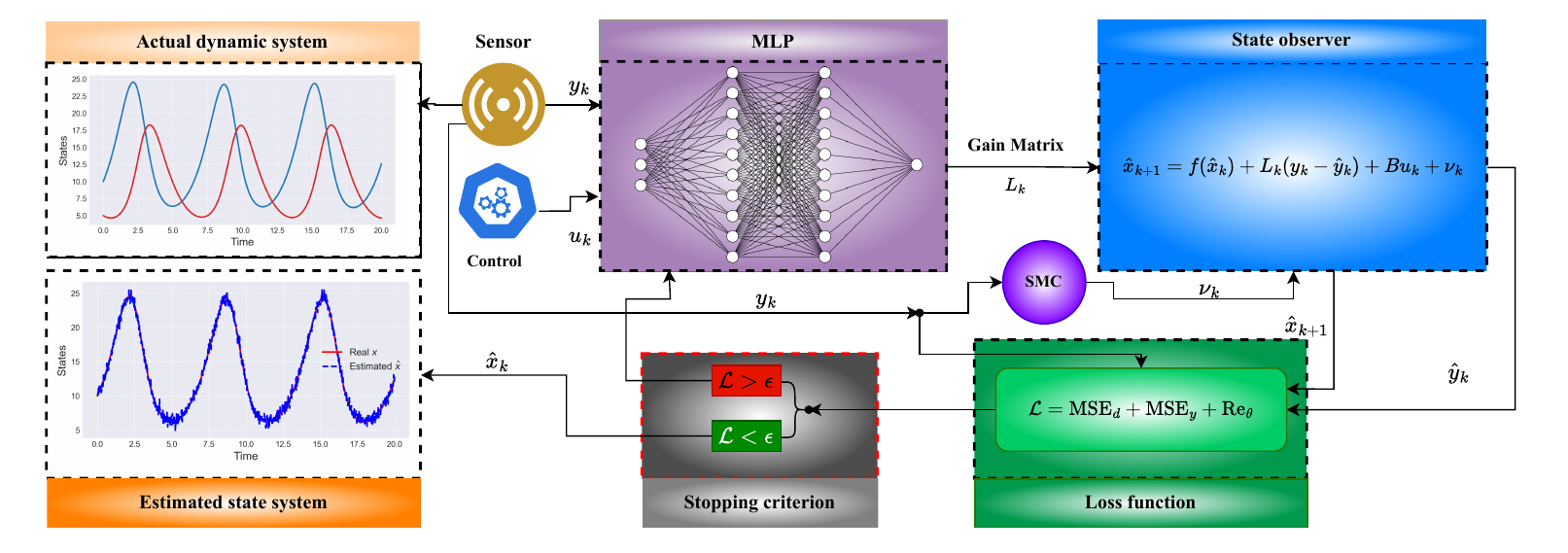}
     \caption{Adaptive nonlinear state observer design-based neural networks. }
     \label{f1}
 \end{figure*}
\subsection{Neural Network}
The time-varying gain matrix \( L_k \) is determined using a neural network that learns the nonlinear relationships between the control input, the system output, and the state errors. At each discrete time step, the neural network processes the current time \( t_k \), control input \( u_k \), and system output \( y_k = h(x_k) \) to generate an adaptive gain matrix \( L_k \), represented as:
\begin{equation}
L_k = \mathcal{N}(t_k, u_k, y_k).   
\end{equation}
Here, the gain matrix \( L_k \) is dynamically updated to account for nonlinearities, disturbances, and noise of the system, improving the observer’s ability to respond quickly and robustly.

The neural network \( \mathcal{N}^{(\ell)} \) is a multilayer feedforward model consisting of \( \ell \) layers. A set of weights and biases in each layer defines the architecture of the network. Let \( \mathbf{W}_i\) and \( \mathbf{b}_i \) represent the weight matrix and the bias vector for the \( i \)-th layer, where \( i = 1, 2, \dots, \ell \). These parameters are collectively denoted as \( \theta_{\ell} := \{\mathbf{W}_i, \mathbf{b}_i\}_{1\leq i \leq \ell }\).

The forward propagation of the network is defined recursively. For an input vector \( z_k = [t_k, u_k, y_k] \), the output of the first layer is given by:
\[
 \mathcal{N}^{(1)} = \sigma(\mathbf{W}_1 z_k + \mathbf{b}_1),
\]
where \( \sigma \) is the activation function applied element-wise. The outputs of subsequent layers are computed recursively as:
\begin{equation}
\mathcal{N}^{(i)} = \sigma(\mathbf{W}_{i} \mathcal{N}^{(i-1)} + \mathbf{b}_{i}), \quad i =  2, \dots, \ell.
\end{equation}

The output of the final layer is the time-varying gain matrix \( L_k \):
\begin{equation}
L_k = \mathcal{N}^{(\ell)}(z_k; \theta_{\ell})
\end{equation}
The activation function \( \sigma \) used in the network is the hyperbolic tangent function, \( \tanh(x) \), which is chosen for its bounded outputs. This helps to mitigate the problem of exploding gradients during training. The derivative of the tanh function is given by:
\[
\frac{d}{dx} \tanh(x) = 1 - \tanh^2(x),
\]
which avoids the vanishing gradient problem by maintaining non-zero gradients during backpropagation.  Furthermore, the tanh function exhibits 1-Lipschitz continuity, which promotes stable learning. This property ensures that small input variations result in bounded output changes, preventing excessive sensitivity and improving the robustness of the training process.
\subsection{Adaptive SMC}
The SMC mechanism operates on the principle of minimizing the output error between the true state and the estimated state. This error is encapsulated in the sliding surface \( s_k \), which is defined as:
\begin{equation}
s_k =  e^{y}_k = h(x_k) - h(\hat{x}_k),
\end{equation}
where \( h(x_k) \) represents the system's true output, and \( h(\hat{x}_k) \) is the output predicted by the state observer. The objective of SMC is to drive this output error \( s_k \) to zero by applying a corrective control action \( \nu_k \), ensuring that the estimation error is minimized despite any external disturbances or uncertainties.

The corrective control action applied by SMC is designed using a hyperbolic tangent function to smooth the control signal and reduce the likelihood of high-frequency oscillations (chattering). The control law is defined as:
\begin{equation}
\nu_k = -\overline{K}_k \tanh(s_k),
\end{equation}
where \( \overline{K}_k \) is an adaptive control gain. The use of \( \tanh(s_k) \) instead of the sign function ensures that the control action is continuous, thus reducing the possibility of chattering, especially during transitions near \( s_k = 0 \). This smooth control action provides a balance between fast convergence and stability.

A key feature of the proposed SMC approach is the use of an adaptive gain \( \overline{K}_k \), which adjusts in real-time based on the magnitude of the sliding surface \( s_k \). The adaptive gain is defined as:
\begin{equation}
\overline{K}_k = K_0 + \alpha \| s_k \|^2,
\end{equation}
where \( K_0 \) is a base gain, and \( \alpha \) is a positive scaling factor. This design increases the control gain when the error is large, leading to a more aggressive error reduction. Conversely, when the error is small, the gain is reduced, minimizing unnecessary control action and preventing over-compensation. Adjusting \( \overline{K}_k \) in response to the magnitude of the error, the control system maintains a high degree of robustness while minimizing control effort.
\subsection{Data and loss functions}
To ensure effective training of the model, the essential dataset must include the following components: time instance \( t_k \), input control torque \( u_k \), sensor measurements \( y_k \) (obtained from real-world data or simulated through ODE/PDE solvers), and initial conditions. This dataset provides the foundation for training the neural network, which aims to minimize the following Mean Squared Error (MSE) loss function:
\begin{equation}
\text{MSE} = \text{MSE}_d + \text{MSE}_y + \text{Re}_{\theta},
\end{equation}
where:
\begin{align}
\text{MSE}_d &= \| \hat{x}_{k+1} - f(\hat{x}_{k}) -Bu_k \|^2, \\
\text{MSE}_y &= \| h(x_k) - h(\hat{x}_{k}) \|^2, \\
\text{Re}_\theta &= \lambda \sum_{i} \| W_i \|^2.    
\end{align}
The MSE loss function incorporates three key components:
\begin{itemize}
\item  Residual Dynamics (\(\text{MSE}_d\)): This term penalizes the solution to satisfy the system dynamics, thereby promoting an accurate estimation of the system's states.
\item  Output Error (\(\text{MSE}_y\)): This component minimizes the discrepancy between the true output and the observed output, ensuring precise tracking of the system's outputs by the observer.
\item  Regularization Term (\(\text{Re}_\theta\)): This term prevents overfitting by imposing a penalty on the weights of the network parameters. This encourages smoother learning of the parameters, enhancing the model's generalization capabilities and robustness.
\end{itemize}
The objective is to solve the following optimization problem:
\[
\min_{\theta_j} \text{MSE}(\theta_j, u_k, y_k, t_k),
\]
where \( \theta_j \) represents the trainable parameters of the neural network. This optimization can be addressed using gradient-based algorithms such as L-BFGS, Adam, and others. These methods iteratively adjust the parameters to minimize the loss function, guiding the observer toward accurate state estimation and improved robustness across various operating conditions.

\subsection{Model training}\label{s5_train}
The proposed neural network can be trained effectively using data from either an ODE-PDE solver or real-world sensors. When using real-world data, measurements can be obtained through physical sensors or by simulating real physical phenomena with models to capture the necessary state data. However, as physical sensors can be costly, minimizing their number while maintaining data quality is critical.

An ODE-PDE solver can be employed to test and validate the model, especially for nonlinear systems. This involves solving the dynamic equations to obtain the output measurements $y$. To initiate this process, an initial condition, $x_0 \in \mathbb{R}^n$, representing the system's initial state, is randomly selected for the nonlinear system. The solver can compute the system's subsequent behavior by applying Euler discretization to this initial condition. Similarly, an initial condition for the observer, $\hat{x}_0 \in \mathbb{R}^n$, should be chosen to differ from $x_0$ to assess the estimator's robustness. Other essential parameters include $K_0$, $\alpha$, $\lambda$, $u$, and $B$.

With these parameters, a training dataset can be constructed as $\mathcal{D}_{\text{train}} = \{ t_k, y_k, u_k \}_{k=1}^{N}$. An additional configuration involves setting the neural network's hyperparameters. Xavier initialization initializes the parameters of the neural network, $\theta$. The Tanh activation function is selected due to its effectiveness in approximating nonlinear functions. For optimization, both the L-BFGS and Adam optimizers are suitable, given their adaptive learning capabilities that help avoid local minima. The number of hidden layers and neurons per layer can be adjusted based on the complexity of the specific problem.

During training, the model with the minimum loss function is saved to obtain the optimal gain matrix, $L_k$. 
\subsection{Model testing}
After determining the optimal training weights, the next crucial step is to evaluate the model's effectiveness. New time instances, denoted as $t_k$, and different initial conditions, where the state of the system is measured and recorded as the test output $y_k$, are used as inputs for this evaluation. Subsequently, the observer gains, $L_k$, essential for estimating the system's state based on its output, are predicted using the obtained training weights. These predicted observer gains are then incorporated into the observer system \eqref{e2}, formulated as a set of differential equations. By numerically solving this system, the estimated states, $\hat{x}_k$, are obtained at different time instances.

The observer's effectiveness is evaluated by calculating the absolute error between the actual and estimated state vectors:
\begin{equation}
\text{Error}_i =| x_{k}^i -\hat{x}_{k}^i|, \quad 1\leq i \leq n
\end{equation} 
To compare the proposed approach with other methods, various metrics are employed, including the Mean Absolute Error (MAE), Root Mean Square Error (RMSE), Mean Squared Error (MSE), and Symmetric Mean Absolute Percentage Error (SMAPE):
\begin{align}
\text{MAE} &= \frac{1}{n} \sum_{i=1}^{n} |x^i - \hat{x}^i| \\   
\text{RMSE} &= \sqrt{\frac{1}{n} \sum_{i=1}^n (\hat{x}^i - x^i)^2} \\
\text{MSE}&=  \frac{1}{n} \sum_{i=1}^{n} (x^i - \hat{x}^i)^2\\
\text{SMAPE}&= \frac{1}{n} \sum_{i=1}^{n} \frac{|x^i - \hat{x}^i|}{(|x^i| + |\hat{x}^i|)/2} \times 100\%
\end{align}
This comparative analysis provides a quantitative measure to assess the observer's accuracy in predicting the system's true state. Lower errors indicate a more precise estimation, while higher errors may necessitate reassessing the model's efficacy.

Through this systematic evaluation, the observer-based model's accuracy in capturing and predicting the dynamics of the underlying system is thoroughly assessed.

\section{Convergence Analysis}
This section presents a rigorous convergence analysis for the proposed adaptive observer framework. The analysis is based on Lyapunov stability theory and establishes sufficient conditions for the exponential convergence of the estimation error to a small neighborhood of the origin. The formal proof relies on the following key assumptions.
\begin{assum}[Bounded Trajectories]
  There exists a compact set \(X \subset \mathbb{R}^{n}\) such that for any solution \(x_k\) of the discrete-time system \eqref{e1} of interest, the state trajectory remains within \(X\) for all \(k \geq 0\):
\[
x_k \in X, \quad \forall k \geq 0.
\]
This assumption ensures that the state trajectories of the discrete-time system are bounded and confined to a finite region of the state space's overall discrete time steps. 
It provides the foundation for the validity of observer design, ensures stability analysis is conducted within a bounded region, and supports robustness to perturbations or modeling inaccuracies.
\label{assum1}
\end{assum}
\begin{assum}[Bounded Jacobians]
    The system function $f$ and observation function $h$ are continuously differentiable on the compact set $\Omega$. Consequently, their Jacobians, $F_x = \frac{\partial f}{\partial x}$ and $H_x = \frac{\partial h}{\partial x}$, are bounded for all $x \in X$:
$$ \left\| \frac{\partial f}{\partial x}\bigg|_x \right\| \leq M_f, \quad \left\| \frac{\partial h}{\partial x}\bigg|_x \right\| \leq M_h $$
for some positive constants $M_f$ and $M_h$ \label{assum2}
\end{assum}
\begin{defn}
   For a nonlinear discrete-time dynamical system~\eqref{e1}, 
   the observability matrix at step \(k\) over a finite horizon \(N\) is defined as:

\[
\text{Obs} =
\begin{pmatrix}
H_k \\
H_{k+1} F_k \\
H_{k+2} F_{k+1} F_k \\
\vdots \\
H_{k+N} F_{k+N-1} \dots F_k
\end{pmatrix}.
\]
The system \eqref{e1} is said to be strongly observable if $\text{rank(Obs)} = n$.
\end{defn}
\begin{assum}[Uniform Observability]
    The system is uniformly observable. This implies that for a finite horizon $N$, there exist positive constants $\alpha, \beta > 0$ such that the observability Gramian formed from the system's Jacobians satisfies:
$$ \alpha I \leq \sum_{i=0}^{N-1} (H_{k+i} \Phi_{k,i})^T (H_{k+i} \Phi_{k,i}) \leq \beta I $$
where  $\Phi_{k,i} = F_{k+i-1} \cdots F_k$ is the state transition matrix of the linearized system and  $F_k = \frac{\partial f}{\partial x}\big|_{\hat{x}_k}$ and $H_k = \frac{\partial h}{\partial x}\big|_{\hat{x}_k}$. This condition is fundamental, as it guarantees that the internal states are reconstructible from the output measurements and ensures the existence of a stabilizing observer gain. \label{assum3}
\end{assum}
\begin{assum}[Neural Network Approximation]
    By the Universal Approximation Theorem, a neural network with sufficient capacity can approximate any continuous function over a compact set. We assume that an ideal, stabilizing observer gain $L_k^*$ exists (guaranteed by Assumption 3). The neural network $N_\theta$ is trained to approximate this gain, resulting in a bounded approximation error $\varepsilon_k$:
$$ \|L_k - L_k^*\| = \|\varepsilon_k\| \leq \varepsilon_{\max} $$
where $\varepsilon_{\max}$ is a small positive constant dependent on the network architecture and training quality.
\label{assum5}
\end{assum}
\begin{thm}[Exponential Stability with NN Approximation Errors]
\label{thm:main_stability}
Under Assumptions~\ref{assum1}--\ref{assum5}  , if the neural network approximation error satisfies $\norm{\varepsilon_k} \leq \varepsilon_0 e^{-\lambda t_k}$ for some $\lambda > 0$, and the adaptive SMC gain satisfies $K_k \geq K_{\min} > 0$, then the estimation error $e_k = x_k - \hat{x}_k$ converges exponentially to a bounded region around zero.
\end{thm}
\begin{prf}
\textbf{Step 1: Error Dynamics Analysis}
The error dynamics are:
\begin{equation}
e_{k+1} = f(x_k) - f(\hat{x}_k) - L_k(h(x_k) - h(\hat{x}_k)) - \nu_k \label{eq:error_dynamics}
\end{equation}

Using the Mean Value Theorem:
\begin{align}
f(x_k) - f(\hat{x}_k) &= F_k e_k + \delta_f(e_k) \label{eq:mvt_f}\\
h(x_k) - h(\hat{x}_k) &= H_k e_k + \delta_h(e_k) \label{eq:mvt_h}
\end{align}

where $\norm{\delta_f(e_k)} \leq \frac{M_f}{2}\norm{e_k}^2$ and $\norm{\delta_h(e_k)} \leq \frac{M_h}{2}\norm{e_k}^2$.

\textbf{Step 2: Decomposition of Gain Matrix}

Decompose the neural network gain as:
\begin{equation}
L_k = L_k^* + \varepsilon_k \label{eq:gain_decomposition}
\end{equation}

where $L_k^*$ is the optimal gain that would achieve desired convergence properties.

\textbf{Step 3: Error Dynamics Reformulation}

Substituting equations~\eqref{eq:mvt_f}, \eqref{eq:mvt_h}, and \eqref{eq:gain_decomposition} into \eqref{eq:error_dynamics}:
\begin{align}
e_{k+1} &= (F_k - L_k^* H_k)e_k - \varepsilon_k H_k e_k + \delta_f(e_k) \nonumber\\
&\quad - L_k^* \delta_h(e_k) - \varepsilon_k \delta_h(e_k) - \nu_k \label{eq:error_reformulated}
\end{align}

\textbf{Step 4: Lyapunov Function Construction}

Consider the Lyapunov function:
\begin{equation}
V_k = e_k^T P_k e_k \label{eq:lyapunov}
\end{equation}

where $P_k \succ 0$ is a positive definite matrix that adapts based on system observability.

\textbf{Step 5: Lyapunov Difference Analysis}

\begin{equation}
\Delta V_k = V_{k+1} - V_k = e_{k+1}^T P_{k+1} e_{k+1} - e_k^T P_k e_k \label{eq:lyapunov_diff}
\end{equation}

Let $A_k = F_k - L_k^* H_k$. Then:
\begin{equation}
e_{k+1} = A_k e_k + w_k \label{eq:error_compact}
\end{equation}

where $w_k$ represents the combined effect of NN approximation errors, higher-order terms, and SMC action:
\begin{align}
w_k &= -\varepsilon_k H_k e_k + \delta_f(e_k) - L_k^* \delta_h(e_k) - \varepsilon_k \delta_h(e_k) - \nu_k \label{eq:disturbance}
\end{align}

\textbf{Step 6: Bound on Lyapunov Difference}

\begin{align}
\Delta V_k &\leq e_k^T (A_k^T P_{k+1} A_k - P_k) e_k + 2e_k^T A_k^T P_{k+1} w_k + w_k^T P_{k+1} w_k \label{eq:lyapunov_bound}
\end{align}

\textbf{Step 7: Design of Optimal Gain $L_k^*$}

Choose $L_k^*$ such that:
\begin{equation}
A_k^T P_{k+1} A_k - P_k \leq -\gamma P_k \label{eq:riccati_condition}
\end{equation}

for some $\gamma > 0$. This can be achieved through the discrete-time Riccati equation approach.

\textbf{Step 8: Bound on Disturbance Term}

For the SMC term, since $\nu_k = -K_k \tanh(s_k) = -K_k \tanh(H_k e_k)$:
\begin{equation}
e_k^T H_k^T \tanh(H_k e_k) \geq \sigma\norm{H_k e_k}^2 \label{eq:tanh_property}
\end{equation}

for some $\sigma > 0$ (property of tanh function).

For the approximation errors:
\begin{align}
\norm{w_k} &\leq \norm{\varepsilon_k} \norm{H_k} \norm{e_k} + \frac{M_f}{2}\norm{e_k}^2 \nonumber\\
&\quad + \norm{L_k^*} \frac{M_h}{2}\norm{e_k}^2 + \norm{\varepsilon_k} \frac{M_h}{2}\norm{e_k}^2 + K_k\norm{H_k e_k} \label{eq:disturbance_bound}
\end{align}

\textbf{Step 9: Sufficient Conditions for Stability}

If the following conditions hold:

\begin{enumerate}
\item \textbf{NN Approximation Quality}: $\norm{\varepsilon_k} \leq \varepsilon_0 e^{-\lambda t_k}$
\item \textbf{SMC Gain Condition}: 
\begin{equation}
K_k \geq K_{\min} > \frac{\norm{\varepsilon_k} \norm{H_k} + \frac{M_f + \norm{L_k^*}M_h + \norm{\varepsilon_k}M_h}{2}\norm{e_k}}{\sigma} \label{eq:smc_condition}
\end{equation}
\item \textbf{Observability Maintenance}: 
\begin{equation}
\gamma > \frac{(\norm{\varepsilon_k} \norm{H_k} + \frac{M_f + \norm{L_k^*}M_h}{2}\norm{e_k})^2}{\lambda_{\min}(P_k)} \label{eq:observability_condition}
\end{equation}
\end{enumerate}

Then:
\begin{equation}
\Delta V_k \leq -\rho V_k + \theta \label{eq:lyapunov_final}
\end{equation}

for some $\rho > 0$ and $\theta \geq 0$ (bounded term).

\textbf{Step 10: Convergence Result}

This implies:
\begin{equation}
V_k \leq e^{-\rho k}V_0 + \frac{\theta}{1-e^{-\rho}} \label{eq:convergence_bound}
\end{equation}

Therefore:
\begin{equation}
\norm{e_k} \leq \sqrt{\frac{\lambda_{\max}(P_k)}{\lambda_{\min}(P_k)}} e^{-\rho k/2}\norm{e_0} + \sqrt{\frac{\theta}{\lambda_{\min}(P_k)(1-e^{-\rho})}} \label{eq:final_bound}
\end{equation}

This establishes exponential convergence to a bounded region.
\end{prf}
\begin{coro}[Ultimate Boundedness]\label{cor:ultimate_bound}
If the neural network approximation error remains bounded ($\norm{\varepsilon_k} \leq \varepsilon_{\max}$), then the estimation error is ultimately bounded with:
\begin{equation}
\limsup_{k \to \infty} \norm{e_k} \leq C \varepsilon_{\max}
\end{equation}
for some constant $C > 0$.
\end{coro}

\section{Applications}
In this section, we extensively assess the proposed method's efficacy in addressing highly complex and nonlinear system dynamics. 

Table~\ref{tab:examples} presents a summary of numerical benchmarks, system complexities, and experimental objectives. To comprehensively evaluate the proposed approach, extensive simulations were performed. Initially, we replicated examples from the state-of-the-art literature, systematically varying parameters to ensure a thorough assessment. Subsequently, as shown in Table~\ref{tab:examples}, we introduced more challenging scenarios by modifying the output $y$ and altering the nonlinearities to examine systems with varying levels of observability complexity. Furthermore, we considered examples with non-differentiable dynamics, which represent significant challenges for existing methods. These complex systems offer valuable insights into the robustness and applicability of the proposed approach.

We employ the following parameter settings for model training and evaluation: $K_0 = 5$, $\alpha = 0.01$, $\lambda = 0.001$. The input of all systems is set to $u = 0$, except for system \eqref{e21}. The neural network architecture contains two hidden layers, each containing $64$ neurons. Other parameters, such as $\delta t$, $T$, $y$, $x_0$, and $\hat{x}_0$, are application-specific and will be explicitly presented in the subsequent sections. 
\begin{table*}[ht!]
\centering
\caption{COMPARISON OF EXAMPLE CASES AND THEIR CHARACTERISTICS.}
\begin{tabular}{llllllll}
\hline
\multicolumn{1}{c}{\multirow{2}{*}{\textbf{Example}}} & \multirow{2}{*}{\textbf{Citation}}                                                  & \multirow{2}{*}{\textbf{Complexity}}                                                                & \multirow{2}{*}{\textbf{Noise ($\omega$)}} & \multirow{2}{*}{\textbf{Input (u)}}        & \multirow{2}{*}{\textbf{Output (y)}}                                  & \multirow{2}{*}{\textbf{Challenge}}        & \multirow{2}{*}{\textbf{Objective}}                                                                                               \\
\multicolumn{1}{c}{}   &    &    &   &    &   &     &        \\ \hline
\multirow{2}{*}{1}                                    & \multirow{2}{*}{\cite{niazi2023learning}}                                                    & \multirow{2}{*}{\begin{tabular}[c]{@{}l@{}}Chaotic system \\ with Gaussian noise\end{tabular}}      & \multirow{2}{*}{\checkmark} & \multirow{2}{*}{\xmark}     & \multirow{2}{*}{$y = x_1$}                                            & \multirow{2}{*}{\xmark}     & \multirow{2}{*}{\begin{tabular}[c]{@{}l@{}}Handles noisy \\ chaos dynamics\end{tabular}}  \\&  &   &  &  &   &   &            \\ \hline
\multirow{2}{*}{2}                                    & \multirow{2}{*}{\cite{sanfelice2015convergence,praly2006new}}                                 & \multirow{2}{*}{\begin{tabular}[c]{@{}l@{}}Invariant state \\ estimation\end{tabular}}              & \multirow{2}{*}{\xmark}     & \multirow{2}{*}{\xmark}     & \multirow{2}{*}{$y = x_1$}                                            & \multirow{2}{*}{\xmark}     & \multirow{2}{*}{\begin{tabular}[c]{@{}l@{}}Works with \\ constant states\end{tabular}}                                       \\  &   &   &  &  &      & &                                   \\ \hline
3                                                     & \cite{miao2023learning}                                                                     & \begin{tabular}[c]{@{}l@{}}Nonlinear system \\ with modified outputs\end{tabular}                   & \xmark                      & \xmark                      & \begin{tabular}[c]{@{}l@{}}$y = x_1$, \\ $y = x_1 + x_2$\end{tabular} & \checkmark                  & \begin{tabular}[c]{@{}l@{}}Modified output \\ handling, superior\\  convergence speed\end{tabular}                                \\ \hline
4                                                     & \cite{sanfelice2011convergence,sanfelice2023convergence}                                                  & \begin{tabular}[c]{@{}l@{}}Nonlinear weakly  \\ observable systems\end{tabular}                 & \xmark                      & \xmark                      & \begin{tabular}[c]{@{}l@{}}$y = x_1$, \\ $y = x_1 + x_2$\end{tabular} & \checkmark                  & \begin{tabular}[c]{@{}l@{}}Handles modified \\ outputs, robust \\ estimation\end{tabular}                                         \\ \hline
5                                                     & Derived from \cite{sanfelice2011convergence}                                                       & Enhanced nonlinearities                                                                             & \xmark                      & \xmark                      & $y = x_1$                                                             & \checkmark                  & \begin{tabular}[c]{@{}l@{}}Tackles higher \\ complexity in system\\ nonlinearity\end{tabular}                                     \\ \hline
\multirow{4}{*}{6}                                    & \multirow{4}{*}{\cite{yu2020liquid,kubalvcik2016predictive,sharma20233,tang2021adaptive}}                                  & \multirow{4}{*}{\begin{tabular}[c]{@{}l@{}}Highly nonlinear \\ dynamic, control input\end{tabular}} & \multirow{4}{*}{\xmark}     & \multirow{4}{*}{\checkmark} & \multirow{4}{*}{$y = x_2$}                                            & \multirow{4}{*}{\checkmark} & \multirow{4}{*}{\begin{tabular}[c]{@{}l@{}}Dynamic constraints, \\ control torque, \\ nondifferentiable \\ dynamics\end{tabular}} \\ &      &    &     &    & &      & \\ &   &   & &   &   &  &    \\ &     &     &     &      &       &      &    \\ \hline
\multirow{3}{*}{7}                                    & \multirow{3}{*}{\cite{ramos2020numerical,peralez2021deep,buisson2023towards,niazi2023learning}} & \multirow{3}{*}{\begin{tabular}[c]{@{}l@{}}Weak differential \\ observability\end{tabular}}         & \multirow{3}{*}{\checkmark} & \multirow{3}{*}{\xmark}     & \multirow{3}{*}{$y = x_1$}                                            & \multirow{3}{*}{\xmark}     & \multirow{3}{*}{\begin{tabular}[c]{@{}l@{}}Robust against noise, \\ compares favorably \\ to other methods\end{tabular}}   \\&      &    &    &  &    &   \\
                         &     &          &     &                                            &        &        &         \\ \hline
\end{tabular}
\label{tab:examples}
\end{table*}
\begin{ex}
Consider the Rössler attractor system, defined by the following system:

\begin{equation}
\left\{ \begin{aligned}
\dot{x}_1 &= -x_2 - x_3, \\
\dot{x}_2 &= x_1 + ax_2, \\
\dot{x}_3 &= b + x_3 (x_1 - c), \\
y &= x_2,
\end{aligned}\right.
\label{e17}
\end{equation}
where the parameters are set to \(a = b = 0.2\) and \(c = 5.7\). Here, $y = Cx$, where $C = [0 \ 1 \ 0]$, thus the desired gain matrix will be a vector in $\mathbb{R}^3$. In addition, a zero mean Gaussian noise term \(w(t) \sim \mathcal{N}(0, 1) \times d\)  with a standard deviation of $(d = 0.01)$ is added to the system. This formulation of the Rössler attractor allows a rich analysis of chaotic behavior influenced by stochastic perturbations.

This system is solved in this study \cite{niazi2023learning}, using a supervised PINN approach. This framework effectively approximates both the linearization transformation and its inverse, facilitating the learning process of the underlying system dynamics.

Here,  our model is trained using a dataset generated from the initial conditions $x_0 = [1, 1, 1]$ and $\hat{x}_0 = [0, 1, 2]$, with a time step $\delta t = 0.001$ over a total simulation time of $T = 10$ s.
 \begin{figure*}[!t]
\centering
\subfloat[Predicted vs. actual state trajectories]{\includegraphics[width=0.44\linewidth]{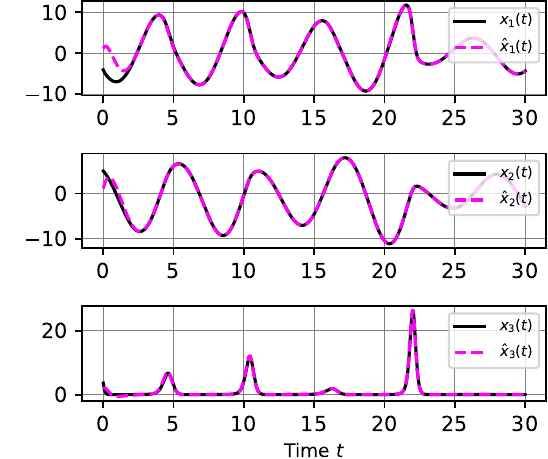}}\hspace{1mm}
\subfloat[Absolute errors over time]{\includegraphics[width=0.45\linewidth]{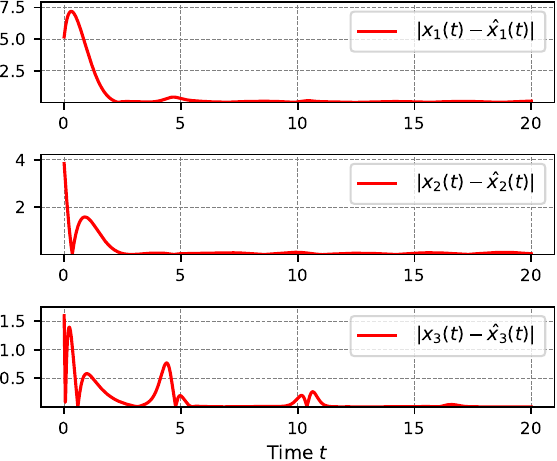}}
\caption{Performance of the proposed method on the Rössler attractor system.  (a) Predicted and actual state trajectories for testing with different initial conditions. (b) Absolute errors over time.}
\label{fig:Rossler}
\end{figure*} 
The best-trained model is obtained for prediction with a minimum loss of $5.8 \times 10^{-3}$. This model is tested with initial conditions $x_0 = [-4, 5, 4]$, significantly different from the initial conditions $\hat{x}_0,  x_0$ used during training, and with an increased time interval  $T=20$s. As illustrated in Fig.~\ref{fig:Rossler}(a), the estimated state graph closely follows the actual state trajectory. Additionally, Fig.~\ref{fig:Rossler}(b) demonstrates the convergence of each state component as time progresses.
These results demonstrate the high performance of our approach in estimating the noisy Rössler attractor system \eqref{e17}. Notably, this approach achieves this accuracy without requiring explicit transformations or linearization, showcasing its robustness to noise and initial conditions and its ability to estimate complex nonlinear dynamics directly. 
\end{ex}
\begin{ex}
The dynamics of the harmonic oscillator \cite{sanfelice2015convergence, praly2006new} are governed by:
    \begin{equation}
\left\{ \begin{aligned}
\dot{x}_1 &= x_2 \\
\dot{x}_2 &= -x_3 x_1 \\
\dot{x}_3 &= 0 \\
y &= x_1
\end{aligned} \label{e18} \right.
\end{equation}
Here, the measured state component is $x_1$, and the goal is to estimate $x_2$ and $x_3$. From the system's dynamics, it is evident that $x_3$ is a constant function. This inherent complexity often challenges state observers, particularly in accurately estimating invariant states. 
The training data is based on the initial conditions $x_0 = [2, -1, 3]$ and $\hat{x}_0 = [1, 1, 2]$. Training is conducted over $T=10$s, achieving a minimum loss of approximately $8.64 \times 10^{-3}$.

\begin{figure*}[!t]
\centering
\subfloat[Predicted vs. actual state trajectories]{\includegraphics[width=0.46\linewidth]{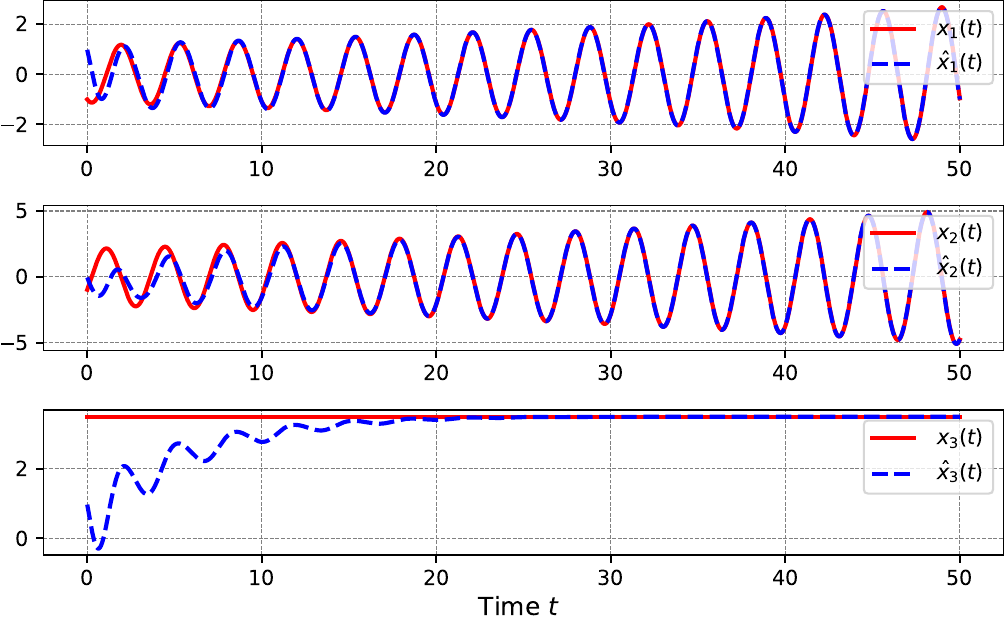}}
\subfloat[Root mean squared error (RMSE)]{\includegraphics[width=0.43\linewidth]{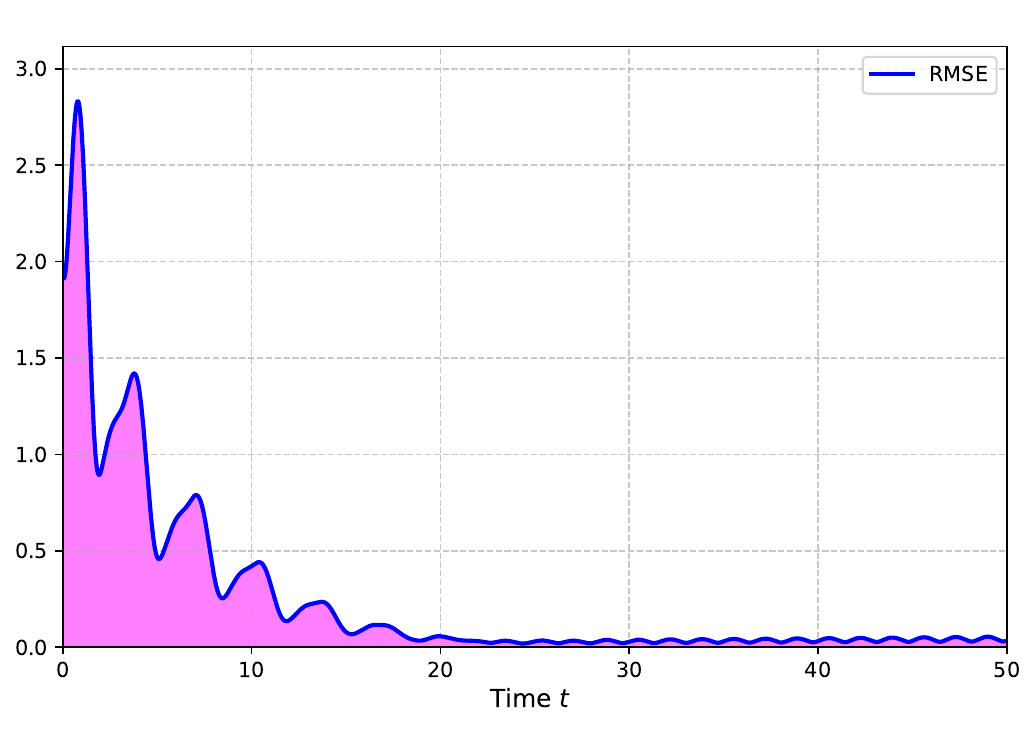}}
\caption{Performance of the proposed method on the harmonic oscillator system. (a) Predicted and actual state trajectories for testing. (b) RMSE over time during testing.}
\label{fig:Harmonic}
\end{figure*}

Fig.~\ref{fig:Harmonic}(a) demonstrates that our model successfully estimates the state of the harmonic oscillator, even for the constant state $x_3$. Fig.~\ref{fig:Harmonic}(b) confirms this convergence by illustrating the RMSE, which tends to zero as time increases. 
\end{ex}
\begin{ex}
To demonstrate the effectiveness of our proposed state estimation model, we consider  another autonomous system, defined by the equations:
\begin{equation}
    \left\{\begin{array}{l}
\dot{x}_1=x_2+\sin x_1 \\
\dot{x}_2=-x_1+\cos x_2\\
y=x_1
\end{array} \right. \label{e19}
\end{equation}
 Previous studies, such as the work \cite{miao2023learning},  have investigated this system within the framework of neural ODEs. In their research, they employed neural ODEs to derive injective mappings and their pseudo-inverses, which facilitated the transformation of the original nonlinear dynamics into a linear representation. 
 
Applying our methodology, the training data is generated from an initial state $x_0 = [1, 1]$ and $\hat{x}_0 = [1, 2]$, and the dynamics are solved over $T = 10$ seconds.
\begin{figure*}[!t]
\centering
\subfloat[Estimated vs. actual states $(y=x_1)$]{\includegraphics[width=0.409\linewidth]{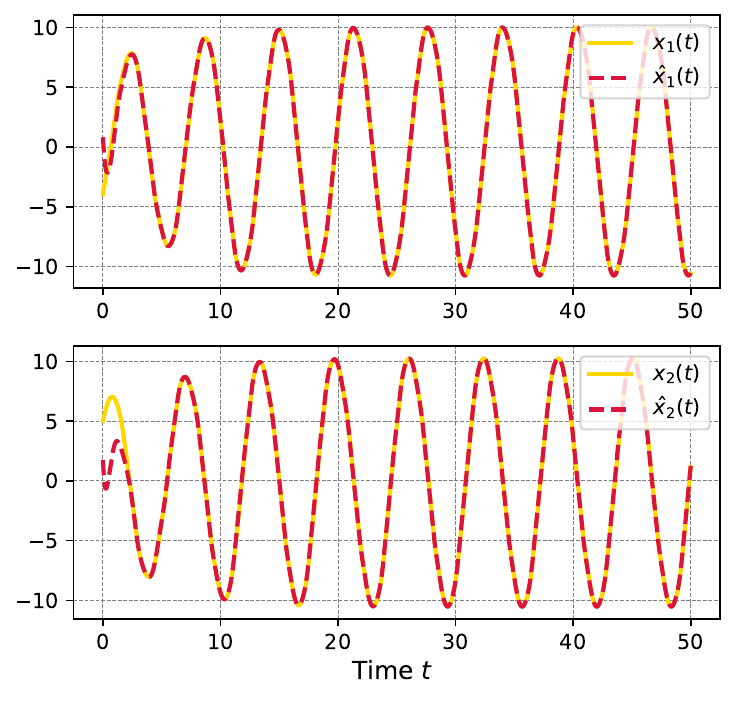}}
\subfloat[RMSE $(y=x_1)$]{\includegraphics[width=0.41\linewidth]{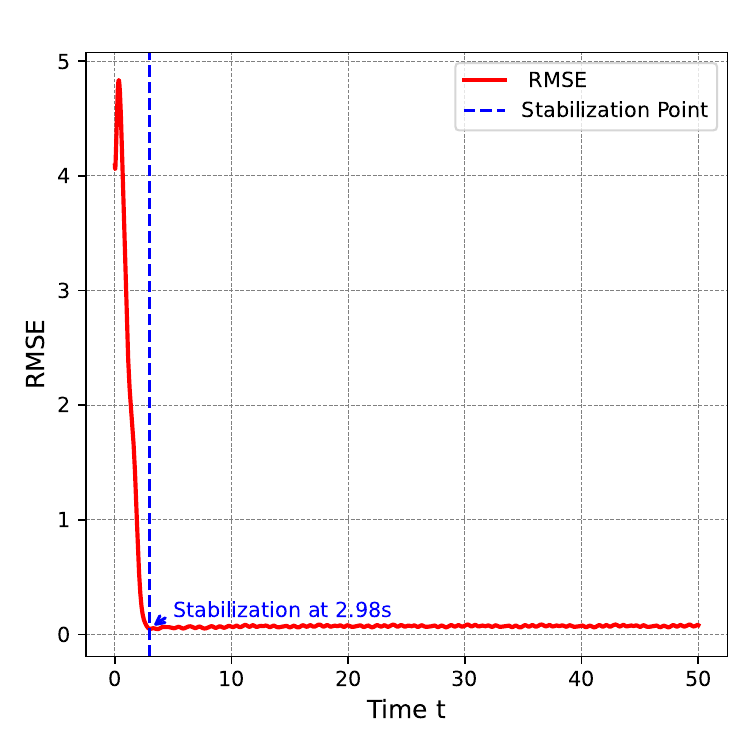}}\\
\subfloat[Estimated vs. actual states $(y=x_1 + x_2)$]{\includegraphics[width=0.409\linewidth]{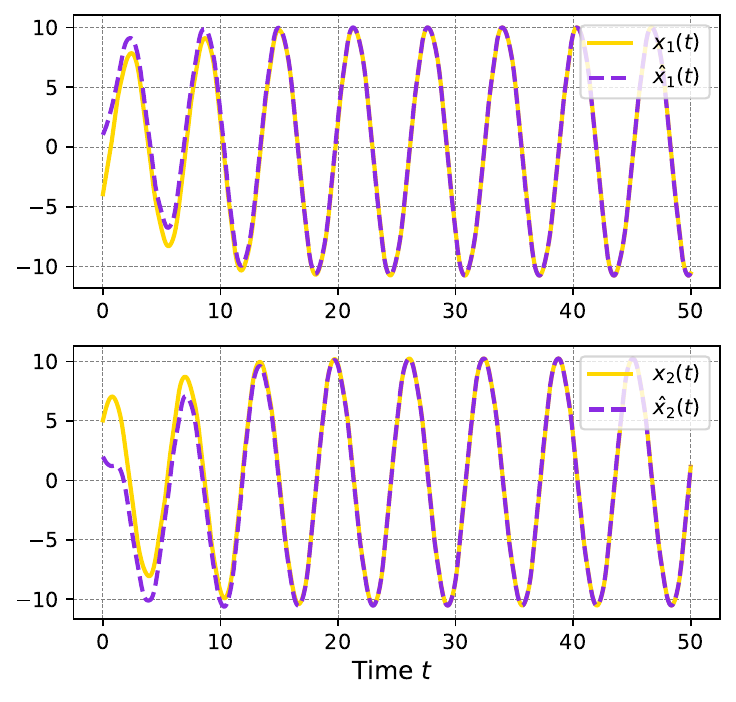}}
\subfloat[Observer trajectories vs actual trajectories $(y=x_1 + x_2)$]{\includegraphics[width=0.41\linewidth]{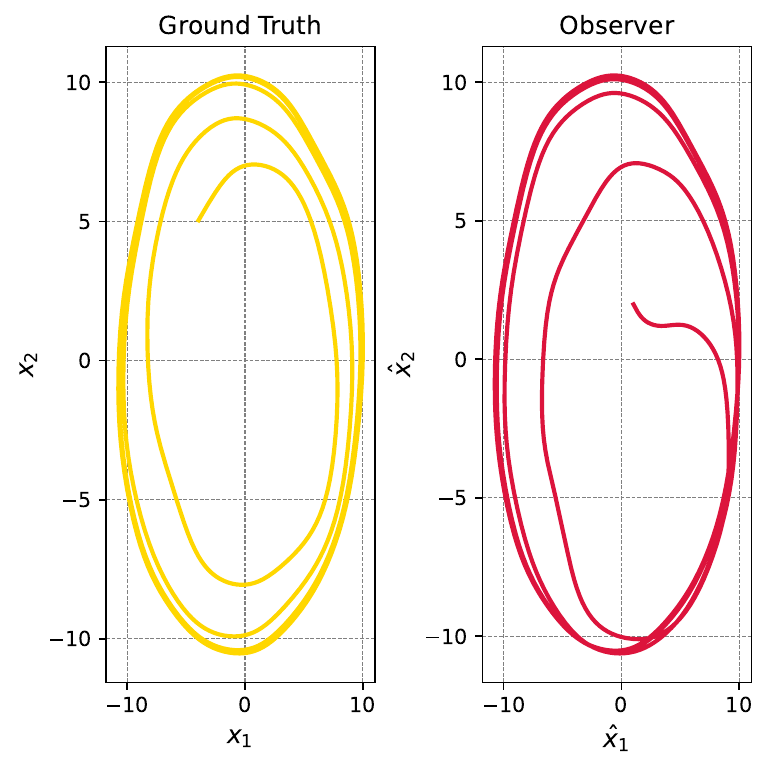}}
\caption{Performance of the proposed method on the autonomous system \eqref{e19}. (a) Estimated and actual state trajectories over time $(y=x_1)$. (b) RMSE convergence $(y=x_1)$. (c) Estimated and actual state trajectories over time with modified output. (d) Observer trajectories vs. actual trajectories with modified output.  }
\label{fig:Autonomous}
\end{figure*}
Fig.~~\ref{fig:Autonomous} demonstrates the accuracy of our proposed model in estimating the state of system \eqref{e19}. The state trajectories converge closely to their true values after a few seconds, confirming the effectiveness of the proposed approach.

Fig.~\ref{fig:Autonomous}(b) illustrates the RMSE and the time of convergence stability. Compared with the stabilization time reported in \cite{miao2023learning}, our model demonstrates superior performance, since it stabilizes the error in $2.98$ seconds (Fig.~\ref{fig:Autonomous}(b)), while the Neural ODEs stabilize within $5$ seconds.  

To increase the complexity of this example, we can modify the output of the system. Instead of measuring $x_1$, we assume that the output measures the sum of both state components:
$$y = x_1 + x_2$$
This modification makes finding a suitable transformation and its inverse a challenging task. Physically, this scenario corresponds to a sensor that measures the sum of $x_1$ and $x_2$, but not the individual components. The performance of our proposed method in this case is demonstrated in Fig.~\ref{fig:Autonomous}(c) and~\ref{fig:Autonomous}(d). These results highlight the efficiency of our model, which does not require any transformation.
\end{ex}
\begin{ex}
Consider the nonlinear academic system presented in \cite{sanfelice2011convergence}:
\begin{equation}
\left\{ \begin{aligned}
 \dot{x}_1 &=x_2 \sqrt{1+x_1^2} \\
\dot{x}_2 &= -\frac{x_1}{\sqrt{1+x_1^2}} x_2^2 \\
y &= x_1
\end{aligned} \right. \label{pra2}
\end{equation}
The work \cite{sanfelice2011convergence,sanfelice2023convergence} establishes observer convergence for nonlinear systems using Riemannian metrics. This approach ensures that the Riemannian distance $d(\hat{x}, x)$ between the system and observer states is non-increasing. Key conditions include the conditional negativity of the Lie derivative $L_fP(x)$, geodesic convexity of the output function level sets, and uniform detectability of the system's linearization. By leveraging the metric's geometry and designing the observer dynamics accordingly, global asymptotic stability of the zero estimation error set can be guaranteed. However, it's worth noting that directly measuring $x_1$ would trivially lead to observability. 
To increase the complexity of this example, we assume  that the output is given by:
$$y = x_1 + x_2$$
In this case, observability becomes non-trivial, and finding the $P$ matrix that satisfies the conditions in \cite{sanfelice2011convergence} may require further investigation.  

The test results of our model are presented in Fig.~\ref{fig:Academic}(a) and~\ref{fig:Academic}(b) for a given initial condition. As confirmed by the trajectories, the $x_1$ trajectory diverges over time. This behavior aligns with the expectation that the system becomes directly observable when measuring $x_1$, as discussed in \cite{sanfelice2011convergence}. However, in our case, we have modified the output measurement. Despite this change, our approach accurately estimates the states, as evidenced in Fig.~\ref{fig:Academic}(a) and~\ref{fig:Academic}(b).
  \begin{figure*}[!t]
\centering
\subfloat[State trajectories]{\includegraphics[width=0.325\linewidth]{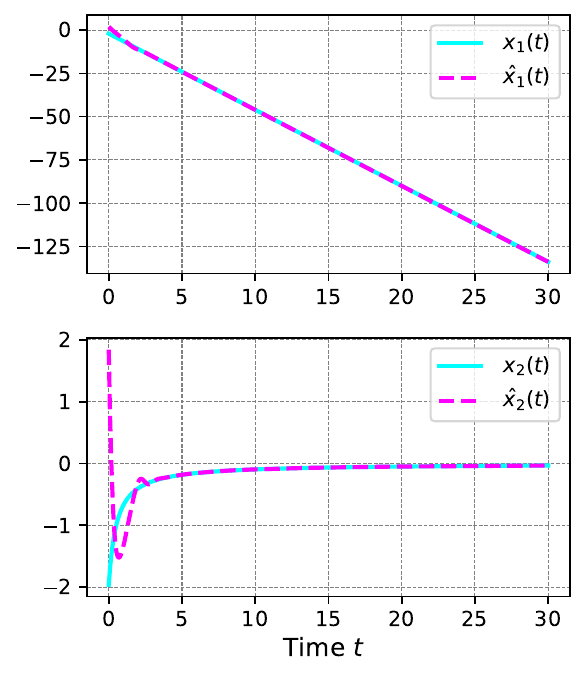}}
\subfloat[Absolute errors]{\includegraphics[width=0.51\linewidth]{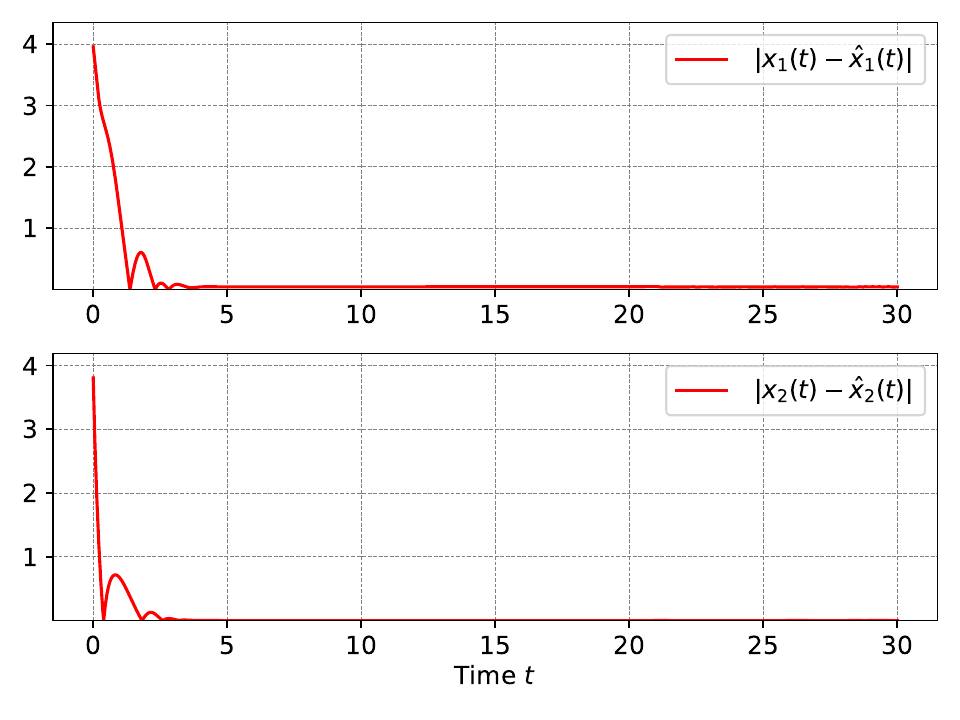}}
\caption{Performance of the proposed method on the academic system \eqref{e21}. (a) State trajectories. (b) Absolute errors over time.}
\label{fig:Academic}
\end{figure*}
\end{ex}
\begin{ex}
Consider another challenging example for linearization-based methods. Consider the nonlinear system:

\begin{equation}
\left\{\begin{aligned}
\dot{x}_1 &= x_2 \sqrt{1+x_2^2} \\
\dot{x}_2 &= -\frac{x_1}{\sqrt{1+x_2^2}} x_2^2 \\
y &= x_1
\end{aligned} \right.\label{pra1}
\end{equation}
This example is derived from system~\eqref{pra2} by modifying the nonlinearity. This modification increases the complexity of the observation of the system.
  \begin{figure*}[!t]
\centering
\subfloat[State trajectories]{\includegraphics[width=0.325\linewidth]{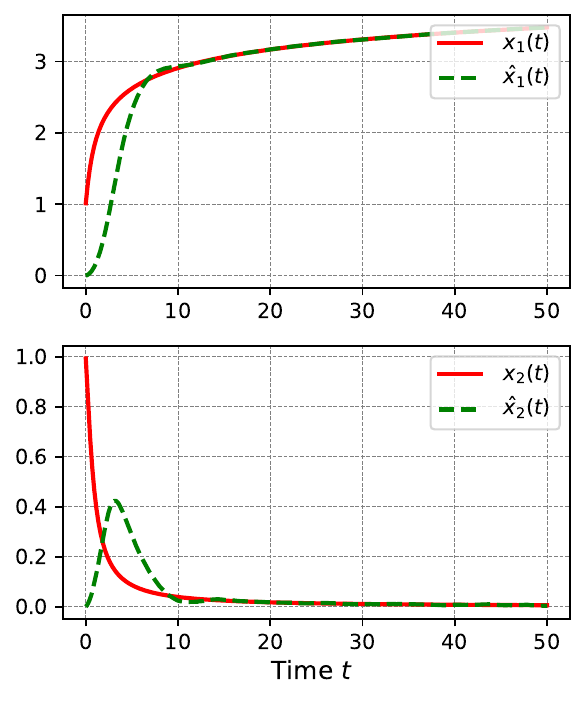}}
\subfloat[Absolute errors]{\includegraphics[width=0.325\linewidth]{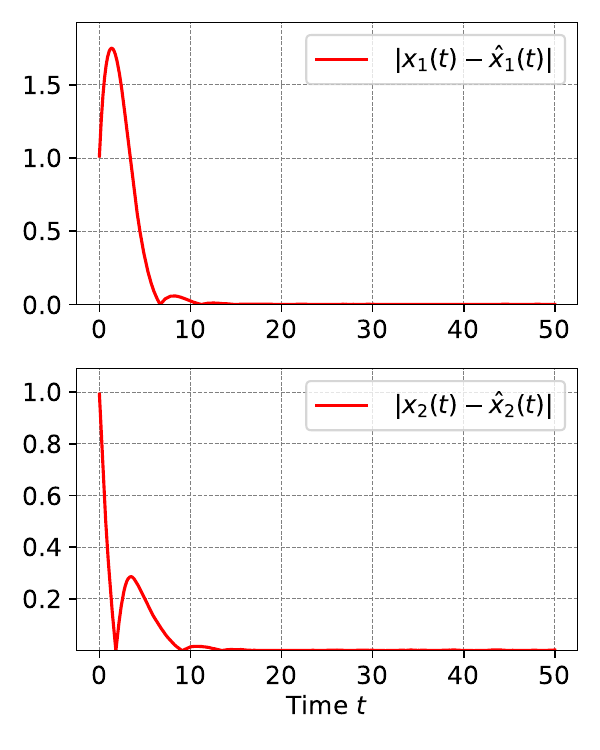}}
\caption{Performance of the proposed method on the academic system \eqref{pra1}. (a) State trajectories. (b) Absolute errors over time.}
\label{fig:Academic1}
\end{figure*}
For a given initial condition, our model estimated the state of this system as shown in Fig.~\ref{fig:Academic1}.
\end{ex}
\begin{ex}
After validating our approach with previous examples, we propose a challenging real nonlinear system benchmark: the Three-Tank System (Fig.~\ref{f_tank}). The following differential equations govern this system:
\begin{equation}
\left\{ \begin{array}{l}
\dot{x}_1 = \frac{u_{1}}{S_T} - \beta_1 \, \text{sign}(x_1 - x_2) \sqrt{|x_1 - x_2|},\\
\dot{x}_2 = \beta_1 \, \text{sign}(x_1 - x_2) \sqrt{|x_1 - x_2|} \\\quad \quad- \beta_2 \, \text{sign}(x_2 - x_3) \sqrt{|x_2 - x_3|},\\
\dot{x}_3 = \frac{u_{2}}{S_T} + \beta_2 \, \text{sign}(x_2 - x_3) \sqrt{|x_2 - x_3|} - \beta_3 \sqrt{x_3},\\
y = x_2,
\end{array} \right.
\label{e21}
\end{equation}
where,
\[
\begin{array}{ll}
x_1, x_2, x_3 & : \text{Water levels in tanks T1, T2, } \\&\text{and T3, respectively,} \\
u_{1}, u_{2} & : \text{Input flow rates from pumps P1 and P2,} \\
S_T & : \text{Cross-sectional area of the tanks,} \\
\beta_1, \beta_2, \beta_3 & : \text{Parameters representing flow characteristics, }\\
&\text{defined as}, 
\beta_i  = \frac{\gamma_{zi} S_p \sqrt{2g}}{S_T}, \quad i \in \{1, 2, 3\},
\end{array}
\]
where $\gamma_{zi}$ are the flow coefficients of the connections between tanks, \( S_p \) is the cross-sectional area of the connecting pipes, and \( g \) is the gravitational acceleration constant.

This system represents a dynamic nonlinear process widely used as a benchmark in control system studies \cite{yu2020liquid,kubalvcik2016predictive,sharma20233,tang2021adaptive,hou2005observing}. The model considers the inflow and outflow dynamics, accounting for the fluid's gravitational and pressure-driven behavior through interconnecting pipes and valves.

The nonlinearity and discontinuity inherent in this system pose significant challenges for the design of observers using conventional methods. These difficulties are primarily due to the root function:
\[Q_{i j}=\beta_i \operatorname{sign}\left(x_i-x_j\right) \sqrt{\left|x_i-x_j\right|}, \quad i,j \in \{1, 2, 3\},\]
 which is nonlinear and non-differentiable at $x_i = x_j$. Additionally, the signum function $\operatorname{sign}(x_i - x_j)$ introduces discontinuities in the system's flow direction at the same point. These characteristics violate the smoothness assumptions required by many observer design techniques, such as KKL or high-gain observers.

While linearization can be employed around a specific operating point, it becomes ineffective near the singularity $x_i = x_j$.

We analyze the Lie derivatives of the output $y$ to assess the system's observability. Considering the measurement of the second tank only, i.e., $y = x_2$, the system is locally observable. However, due to the coupling and nonlinearities, the observability rank may degenerate under specific operating conditions, particularly when $x_i = x_j$.

These measurements can be acquired using physical sensors, such as the MPX 2010 from Motorola \cite{ahmed2020three}.
\begin{figure}
\centering
    \includegraphics[width=1\linewidth]{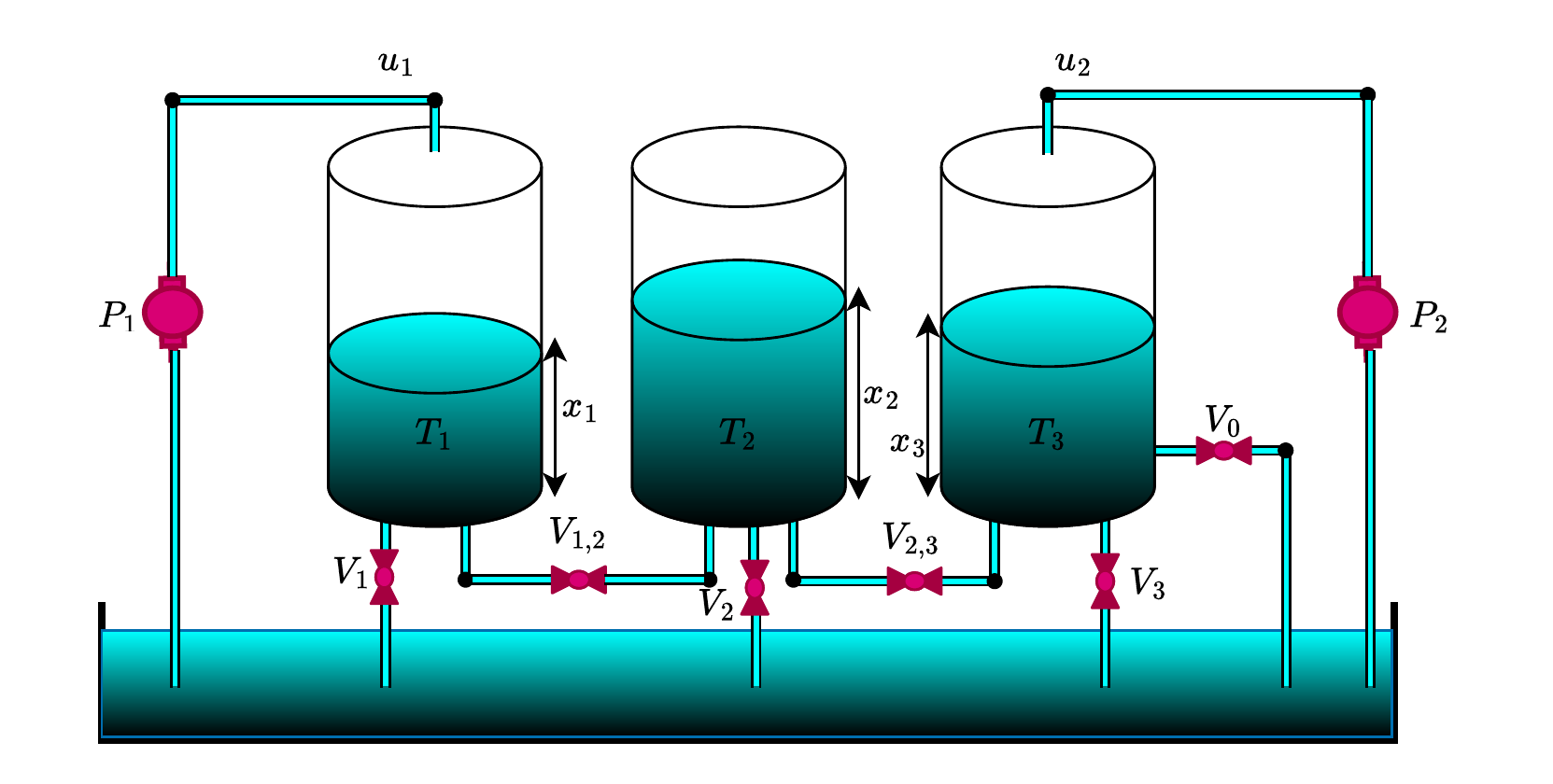}
    \caption{Three tank  system \eqref{e21}.}
    \label{f_tank}
\end{figure}

We conducted two main experiments to evaluate the performance of our model for the three-tank system. Table~\ref{t1} summarizes each experiment's control torque profiles.
\begin{table}
\centering
\caption{CONTROL EXPERIMENT TASKS FOR THE THREE-TANK SYSTEM.}
\begin{tabular}{l|l}
\hline
\textbf{Experiment} & \textbf{Control Torque ($u_1$, $u_2$)} \\ \hline
(1) & $u_1 = u_2 = 0$ \\ \hline
(2) & 
$u_1 = u_2 =   
\begin{cases} 
u_{\text{max}}, & \text{if} \sin(5 \pi f t) > 0, \\
u_{\text{min}}, & \text{otherwise}.
\end{cases}$
 \\ \hline
\end{tabular}

\label{t1}
\end{table}
\begin{itemize}
    \item Experiment (1): In the first experiment, no external control inputs were applied to the system. This means both control torques were set to zero, denoted as $u_1$ and $u_2$. Consequently, the tanks only experienced the natural flow dynamics between them. The resulting state trajectories and absolute errors are shown in Fig.~\ref{fig:ThreeTank}(a) and~\ref{fig:ThreeTank}(b), respectively. As expected, the tank levels converge to zero due to the absence of inflows. Our proposed approach demonstrates good performance in this scenario.
    \item  Experiment (2): The second experiment introduced a control input, denoted as $u$, implemented as a square wave signal. This signal alternated between two values, $u_{max}$ and $u_{min}$, with the switching governed by the sine function $sin(5\pi ft)$. This periodic behavior necessitated adding a non-negativity constraint to both the actual and observed states:
$$x_i \geq 0, \quad \hat{x}_i \geq 0, \quad \forall i = 1, 2, 3$$
The resulting state trajectories and absolute errors are presented in Fig.~\ref{fig:ThreeTank}(c) and~\ref{fig:ThreeTank}(d), respectively. Despite the additional constraint and the dynamic control input, our estimator successfully approximated the state trajectories with high accuracy.
\end{itemize}
\begin{figure*}[!t]
\centering
\subfloat[State trajectories for Experiment (1) (No control)]{\includegraphics[width=0.31\linewidth]{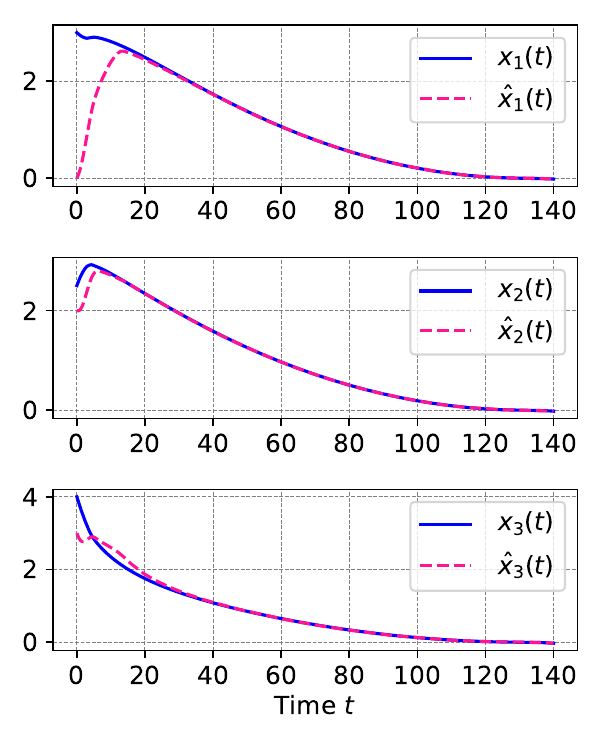}}
\subfloat[Absolute errors for Experiment (1)]{\includegraphics[width=0.335\linewidth]{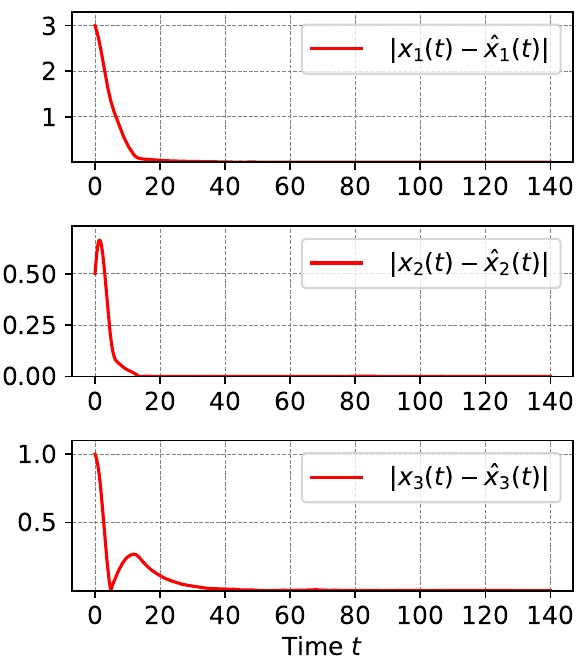}}\\
\subfloat[State trajectories for Experiment (2) (With control)]{\includegraphics[width=0.33\linewidth]{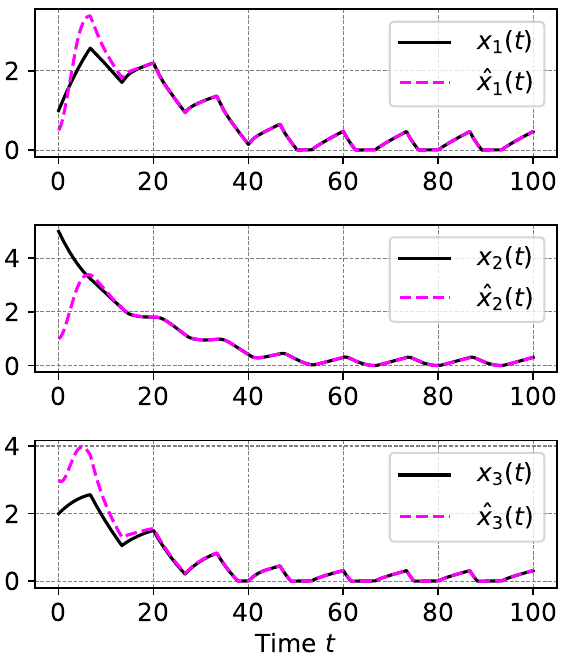}}
\subfloat[Absolute errors for Experiment (2)]{\includegraphics[width=0.33\linewidth]{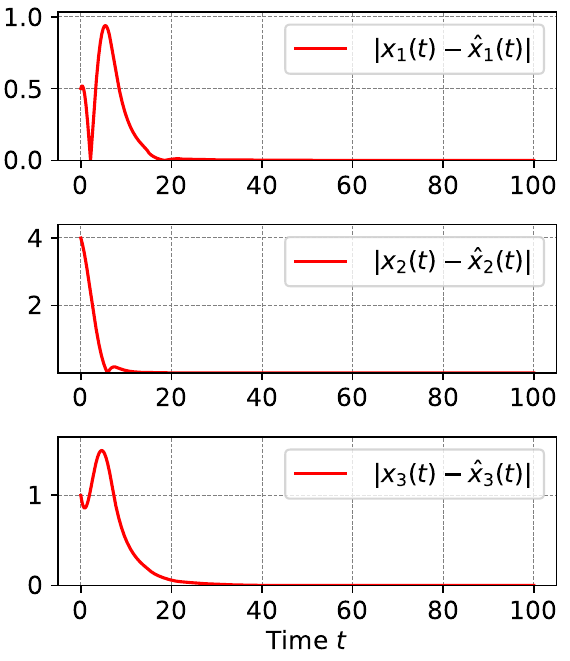}}
\caption{Performance of the proposed method on the Three-Tank system. (a) and (b) correspond to Experiment (1), where no control input is applied. (c) and (d) correspond to Experiment (2), where a square wave control input is applied. The results demonstrate the model's ability to accurately estimate the system's state under both scenarios.}
\label{fig:ThreeTank}
\end{figure*}
\end{ex}
\begin{ex}\label{ex3}
The reverse Duffing oscillator:
\begin{equation}
\left\{\begin{array}{l}
\dot{x}_1=x_2^3 \\
\dot{x}_2=-x_1\\
 y=x_1
\end{array} \right.
\end{equation}
The solutions of this system are known to evolve within invariant compact sets. This system presents challenges in state estimation due to its weak differential observability. This property makes it difficult to accurately reconstruct the system's state from the available output measurements \cite{ramos2020numerical,peralez2021deep,buisson2023towards,niazi2023learning}.

In this study, we evaluated the performance of our approach against several learning-based methods under two scenarios: (1) systems without noise and (2) systems subjected to a noise term \( w(t) \sim \mathcal{N}(0, 0.01) \). This allowed us to assess the robustness of each observer.  

To validate the effectiveness of our method, we compared it with three state-of-the-art approaches: Supervised Neural Networks (Supervised NN), Unsupervised Autoencoders (Unsupervised AE), and Supervised Physics-Informed Neural Networks (Supervised PINNs). The training process for these methods typically relies on a set of predefined training trajectories. However, for our approach, we utilized a set of initial conditions as a training dataset to learn an optimal gain that ensures convergence. For consistency and fairness, we adopted the same neural network architecture described in \cite{niazi2023learning}. 

To measure the efficiency of these approaches, various error metrics such as MSE, RMSE, MAE, and SMAPE were evaluated. The results are presented in Table~\ref{tab:error_comparison}. 
\begin{table*}[ht!]
\centering
\caption{COMPARISON OF ERROR METRICS FOR DIFFERENT MODELS WITH AND WITHOUT NOISE. THE TRAINING DATA FOR OUR APPROACH ARE DERIVED FROM A SET OF INITIAL CONDITIONS, SPECIFICALLY $(x_0, \hat{x}_0) \in [-1, 1] \times [-2, 2]$.}
\begin{tabular}{@{}lcccc@{}}
\toprule
\textbf{Error}&\textbf{Supervised NN}&\textbf{Unsupervised AE}&\textbf{Supervised PINN} & \textbf{Our approach} \\
\midrule
\multicolumn{5}{c}{\textbf{Without Noise}} \\
\midrule
MSE  & $5.972 \times 10^{-5}$ & $4.467 \times 10^{-5}$ & $7.064 \times 10^{-5}$ & $\bm{2.226 \times 10^{-6}}$ \\
RMSE & $7.728 \times 10^{-3}$ & $6.683 \times 10^{-3}$ & $8.405 \times 10^{-3}$ & $\bm{1.492 \times 10^{-3}}$ \\
MAE  & $5.380 \times 10^{-3}$ & $4.769 \times 10^{-3}$ & $5.868 \times 10^{-3}$ & $\bm{8.884 \times 10^{-4}}$ \\
SMAPE \% & $3.474$ & $2.956$ & $3.839$ & $\bm{0.677}$ \\
\midrule
\multicolumn{5}{c}{\textbf{With Noise $\omega \sim \mathcal{N}(0, 0.01)$}} \\
\midrule
MSE  & $7.452 \times 10^{-5}$ & $6.123 \times 10^{-5}$ & $8.934 \times 10^{-5}$ & $\bm{3.521 \times 10^{-6}}$ \\
RMSE & $8.634 \times 10^{-3}$ & $7.824 \times 10^{-3}$ & $9.450 \times 10^{-3}$ & $\bm{1.875 \times 10^{-3}}$ \\
MAE  & $6.124 \times 10^{-3}$ & $5.883 \times 10^{-3}$ & $6.724 \times 10^{-3}$ & $\bm{1.025 \times 10^{-3}}$ \\
SMAPE \% & $4.235$ & $3.876$ & $4.562$ & $\bm{0.852}$ \\
\bottomrule
\end{tabular}
\label{tab:error_comparison}
\end{table*}
  The results confirm the high accuracy of our method compared to other approaches.
\end{ex}
\section{Conclusion}
This paper introduced a nearly universal approach to the design of software sensors for nonlinear dynamical systems, combining neural networks with adaptive SMC. The proposed approach demonstrated several advantages over existing methods, including:
\begin{enumerate}
    \item Noise Robustness: The model is robust to noise, ensuring accurate estimations even under challenging conditions.
\item No Transformation or Linearization Required: Unlike many traditional methods, this approach does not rely on transformations or linearization, making it suitable for complex nonlinear systems.
\item Fast Training and Convergence: The adaptive SMC with online gain adjustment enables rapid error stabilization and efficient training in discrete time, facilitating real-time implementation and online learning.
\item General Applicability: With minimal restrictive conditions, this method is highly adaptable and applicable to a wide range of nonlinear dynamical systems.
\item High Performance: The proposed approach consistently outperforms existing accuracy, stability, and adaptability methods.
\end{enumerate}
Through rigorous testing on benchmark systems, we validated the robustness and effectiveness of this framework. It surpasses the existing methods and introduces a scalable and reliable solution for nonlinear state estimation.

Despite promising results, further work should explore the performance of the proposed approach under various noise profiles, ensuring its robustness under diverse conditions. Additionally, formal proofs for the stability of the adaptive SMC and the convergence of the methodology in systems with nondifferentiable dynamics would provide critical theoretical insights. Extending this work to higher-dimensional and multivariable systems would broaden its applicability to large-scale real-world scenarios.

This work contributes valuable insights into the design and implementation of effective state estimators for complex systems, paving the way for further advancements in this critical field of research.

\bibliographystyle{IEEEtran}

\vfill
\end{document}